\documentclass[12pt]{article}
\usepackage{amsmath,amsfonts,amssymb,amsthm,amscd}
\usepackage{a4wide}
\title{\it {Special points and Poincar\'e bi-extensions }}
\author{by Daniel Bertrand, with an {\it Appendix} by Bas Edixhoven}
\newtheorem{Theorem}{Theorem}

\newcommand{\C}{\mathbb C}

\newcommand{\Q}{\mathbb Q}
\newcommand{\Z}{\mathbb Z}

\newcommand{\G}{\mathbb G}

\newcommand{\Pic}{\mathrm{Pic}}
\newcommand{\calO}{\mathcal{O}}
\newcommand{\divisor}{\mathrm{div}}
\newcommand{\Div}{\mathrm{Div}}
\newcommand{\quot}{\mathrm{quot}}
\newcommand{\End}{\mathrm{End}}
\newcommand{\Norm}{\mathrm{Norm}}

\date{April 2011}
\begin{document}
\maketitle

 \noindent
 The context is the following :   
 
 \medskip

i) in a joint project with D. Masser,  A. Pillay and U. Zannier \cite{BMPZ}, we aim at extending to semi-abelian schemes  the Masser-Zannier approach \cite{MZ} to Conjecture 6.2 of R. Pink's  preprint  \cite{Pk}; this conjecture also  goes under the name ``\,{\it Relative Manin-Mumford}\,".   Inspired by Anand  Pillay's suggestion that the semi-constant extensions of \cite{BP} may bring trouble, I found a counter-example, which is described   in Section 1 below. 

\smallskip
ii) at a meeting in Pisa end of March, Bas  Edixhoven found a more concrete way of presenting  the counter-example,  with the additional advantage that the order of the involved torsion points can be controlled in a precise way : this is the topic of the Appendix. 

\smallskip
iii) finally, I realized that when rephrased in the context of mixed Shimura varieties, the construction, far from providing a counter-example, actually {\it supports} Pink's  general Conjecture 1.3 of \cite{Pk}; a sketch of this view-point is given in Section 2.

\section{A counter-example to relative Manin-Mumford ...}

\medskip
This counterexample is provided by a ``Ribet section" on a semi-abelian scheme $B/X$ of relative dimension 2 over a base curve $X$. Roughly speaking,   given  an elliptic curve $E_0/\C$  {\it with complex multiplications} and using an idea of L. Breen, K. Ribet  constructed a non-torsion point $\beta_0$ with strange divisibility properties on  any given non isotrivial extension $B_0$ of $E_0$ by $\G_m$, cf. \cite{JR}. In the relative situation $B/X$, the very same construction  yields\,:

\begin{Theorem} Let $B/X$ be  a non constant (hence non isotrivial) extension of $E_0 \times X$
 by $\G_m$. There exists  a section $\beta : X \rightarrow B$ which does not factor through any proper closed subgroup scheme of $B/X$, but whose image $Y := \beta(X)$ meets the torsion points of the various fibers  of $B/X$ infinitely often (so, Zariski-densely, since $X$ is a curve).
 \end{Theorem}

More precisely, let $X$ be a smooth connected affine curve defined over (say) $\C$, with function field $K := \C(X)$. We may have to delete some points of $X$, or consider finite covers of $X$, but will still denote by $X$ the resulting curve. We write $x$ for the generic point of $X$,  i.e. $K = \C(x)$, and $\xi  \in X(\C)$ for its closed points.  We start  with  a CM elliptic curve $E_0/\C$, denote by $\hat E_0  \simeq E_0$ its dual, and fix an antisymmetric isogeny
$$\varphi : \hat E_0 \rightarrow E_0.$$
This means that its transpose $\hat \varphi   : \hat E_0 \rightarrow E_0$ is equal to $- \varphi$, i.e., in the identification $\hat E_0 \simeq E_0$, that $\varphi$ is a totally imaginary complex multiplication. 

\smallskip
We consider the constant elliptic schemes $E = E_0 \times_\C X, \hat E = \hat E_0 \times_\C X$, and fix a non constant section $q : X \rightarrow \hat E$. In particular, $q$ is not a torsion section, and the semi-abelian scheme $B/X$ attached to $q$ is a non constant\footnote{~In particular, $B/X$ is a semi-constant semi-abelian variety in the sense of  \cite{BP}. However, the counter-example to Lindemann-Weierstrass given there is of a different nature; see also Remark 1.(ii) below.},   hence non isotrivial, extension  of   $E/X$ by $\G_m$; conversely, any  such $B$ is of this type. Let $\pi : B \rightarrow E$ be the corresponding $X$-morphism. Extending the construction of  \cite{JR}, we will attach to these data $(q, \varphi)$ a section $\beta $ of $B/X$ such that the section $p := \pi \circ \beta$ of $E/X$ satisfies $p= 2 \; \varphi \circ q$.  Furthermore, $\beta$ will have the following ``{\it lifting property}\;" :  {\it  for any $\xi \in X(\C)$ such that $p(\xi)$ is a torsion point on the fiber $E_{\xi} \simeq E_0$ of $E \rightarrow X$,  its lift $\beta(\xi)$ is automatically a torsion point on the fiber $B_\xi$ of $B \rightarrow X$}.

\medskip
Before describing this construction, let us show that such a section $\beta$ does satisfy the conditions  announced in  Theorem 1 :

\smallskip
\noindent
$\bullet$ on  the one hand, since $q \in \hat E(X)$ has infinite order, the only  proper closed subgroup schemes of $B$ projecting onto $X$ are contained in finite unions of translates of $\G_m \times_\C X$ 
and of fibers of the projection; and the section $p = 2 \; \varphi \circ q \in E(X)$ too has infinite order. Therefore, $\beta$, which projects  to $p$ on $E$, cannot factor through any proper closed subgroup scheme of $B/X$.

\smallskip
\noindent
$\bullet$ on the other hand, since $q$ is not a constant section of $\hat E/X$, neither is $p \in E(X)$, and the set $X^E_{tor} = \{\xi \in X(\C), p(\xi)$ is a torsion point on $E_\xi  \simeq E_0\}$ is infinite. But  by the lifting property, this set coincide the set $X^B_{tor} = \{\xi \in X(\C), \beta(\xi)$ is a torsion point on $B_\xi \}$. Therefore,  the curve $Y = \beta(X) \subset B $  meets the set of torsion points of the various fibers  of $B \rightarrow X$ Zariski-densely.

\bigskip
To perform the construction of $\beta$, we go to the generic fiber $B_x := B \otimes_S K := B_K$ of $B/S$, consider the non-constant point $q(x) := q_K \in \hat E_K = \hat E_0 \otimes_\C K$, and recall the construction of the Ribet point $\beta_K \in B_K(K)$ attached to $q_K$ and to the antisymmetric isogeny $\varphi : \hat E_K \rightarrow E_K$, with transpose $\hat \varphi = - \varphi$. There are two ways to describe $\beta_K$  :

\smallskip
(i) The first one \cite{JR} goes as follows : consider the pullback  
$$\varphi^*B_K \in Ext(\hat E_K, \G_m) ~{\rm of}Ê~ B_K \in Ext(E_K, \G_m)$$
under $\varphi$, and denote  again by $\varphi : \varphi^*B_K \rightarrow B_K$ the natural  extension of $\varphi$ to $\varphi^*B_K$. Since $B_K$ is parametrized by $q_K$, $\varphi^*B_K$ is parametrized by the point $\hat \varphi(q_K)$ of the dual $E_K$ of $\hat E_K$. Now, choose an arbitrary $K$-rational point $t_K$ in the fiber above $q_K$  of the extension $\varphi^*B_K$; in particular, its image  $t^1_K := \varphi(t_K) \in B_K(K)$ satisfies $\pi_K(t^1_K) = \varphi(q_K)$, where $\pi_K = \pi(x) : B_K \rightarrow E_K$. The  point $t_K$ defines a one-motive $M_K: \Z \rightarrow \varphi^*B_K$, whose Cartier dual $\hat M_K : \Z \rightarrow B_K$ is given by a point $t^2_K \in B_K(K)$ projecting to $\hat \varphi (q_K) \in E(K)$. Finally, set $\beta_K = t^1_K - t^2_K$ : this point of $B_K(K)$ is {\it independent of the   choice of the auxiliary point $t_K$} above $q_K$ (this is clear on the symmetric definition of duals given in Hodge III), and its image under $\pi_K$ is the point $  \varphi (q_K) - \hat \varphi (q_K) = 2\varphi(q_K) := p_K$ of $E_K(K)$.

\smallskip
(ii) The second one  \cite{Bd}  is more geometric (and will actually not be used here): consider the Poincar\'e bundle ${\cal P}_K$ on $\hat E_K \times  E_K$, 
rigidified above $(0, 0)$. Since $\varphi$ is antisymmetric, the square of its restriction to the graph $\Phi_K \subset \hat E_K \times  E_K$ of $\varphi$   is trivial. Up to a 2- isogeny, we therefore get a unique non-zero $K$-regular section $\sigma_\varphi : {\Phi}_K \rightarrow {\cal P}_K |_{\Phi_K}$. Now $B_K$ (plus a zero section)  identifies with  the restriction of ${\cal P}_K$ to $ \{q_K \} \times E_K \simeq E_K$, and its fiber over  $p'_K  := \varphi(q_K)$ with the fiber of ${\cal P}_K$ over $(q_K, p'_K) \in \Phi_K(K)$.  We then set $\beta'_K  = \sigma_\varphi\big((q_K, p'_K)\big) $, and  view $\beta'_K$  as a point of $B_K(K)$ above $p'_K$. Up to multiplication  by 2, this is the same point as the $\beta_K$ above.

\medskip
Restricting $X$ if necessary, we can extend this point $\beta_K$  to a  section $\beta  : X  \rightarrow B$, which may be called   the Ribet section  attached to $\varphi$ of the   semi-abelian  scheme $B/X$  defined by the non-constant section $q$ of $ \hat E/X$ we  had started with.  More precisely, we can extend the auxiliary point $t_K$ of the first construction to a section $t$ of $\varphi^*B/X$, and repeat the whole process  over $X$ (minus some points), getting in particular a smooth one-motive $M/X$, sections $t^1, t^2$ over $X$, etc. By definition, $\pi \circ \beta$ is the section $p = 2 \; \varphi \circ q$ of $E/X$ extending $p_K$ over $X$, and it remains to show that $\beta$ satisfies the ``lifting property".

\medskip
So, let $\xi \in X(\C)$ be a point such that $p(\xi)$ is a torsion point on the fiber $E_\xi$ ($\simeq E_0$) of $E \rightarrow X$. We must show that $\beta(\xi)$ is a torsion point on the fiber $B_\xi$ of $B \rightarrow X$. By the relation $p(\xi) = 2 \varphi(q(\xi))$, $q(\xi)$  too is a torsion point on $\hat E_\xi$ (in passing, this shows that
 $B_\xi$ is an isotrivial extension).  So, among the points which lie on  the fiber $(\varphi^*B)_\xi$ of  $\varphi^*B/X$ above $\xi$  (which is the pull-back
 $$  \varphi^* B_\xi  \in Ext(\hat E_\xi, \G_m) ~{\rm of}Ê~ B_\xi \in Ext(E_\xi, \G_m)~Ê),$$ 
 and which  project  to $q(\xi) \in (\hat E_\xi)_{tor}$,  we now have not only the value $t(\xi)$ of the section $t : X \rightarrow \varphi^*B$ at $\xi$, but also plenty of torsion points of  the  complex semi-abelian variety  $\varphi^* B_\xi $.   Choose one of them, and call it  $\tilde t_\xi$. Since $\tilde t_\xi$ and $t(\xi)$ differ by an element of $\G_m$, the first construction, whether applied to $t(\xi)$ or to $\tilde t_\xi$, will yield the {\it same} point $\beta_\xi \in B_\xi$, with $\pi(\beta_\xi) = p(\xi)$. Using $t(\xi)$, we see that $\beta_\xi = \beta(\xi)$;  using the torsion point $\tilde t_\xi$ and the fact that $\hat \varphi (q(\xi))$ is a torsion point, we see that the weight filtrations of the corresponding complex one-motive $\tilde M_\xi$, hence of its dual $\hat {\tilde M}_\xi$, split  up completely up to isogeny. Consequently,   the   points ${\tilde t}^1_\xi, \; {\tilde t}^2_\xi$ and $\beta_\xi$ associated to $\tilde t_\xi$ by the first construction  are all torsion points, and  $\beta(\xi)$ is indeed a torsion point of $B_\xi$.

\bigskip
\noindent
{\bf Remark 1} : i) ({\it from $X$ to $ \hat E_0$}) : let ${\cal B}$ be the ``universal" extension of $E_0$ by $\G_m$, viewed as a group scheme over  $\underline{Ext}^1(E_0, \G_m) \simeq   \Pic^0_{E_0/\C} = \hat E_0$. The extension $B$ attached to the section $q : X \rightarrow \hat E_0$  is  the pull-back of ${\cal B}$  under $q$. Choosing $X = \hat E_0$, and $q = $ the identity map, so that $K = \C(\hat E_0)$, we can therefore restrict to the case where   $K = \C(\hat E_0)$ and $q_K$ is  the generic point of $\hat E_0(K)$. The Appendix - and most of \S 2 - concerns  this generic case $$X = \hat E_0, q = id, B = {\cal B}.$$

\smallskip
\qquad ii)  {\it (when $g > 1$) }: Ribet sections $\beta$ can be defined over any abelian scheme $A/X$, of relative dimension $g$, which admits an antisymmetric isogeny $\varphi : \hat A \rightarrow A$. If $g = 1$, this forces $A$ be to be iso-constant, hence the $E_0/\C$ above. But as soon as $g > 1$, there are examples of simple non constant $A/X$ with such a  $\varphi = - \hat \varphi $.  The section $\beta$ attached to $\varphi$ and to a section $q \in \hat A(X)$ will again satisfy the ``lifting property". However, in order to ensure that  the set $X^A_{tor} = \{\xi \in X(\C), \pi \circ \beta (\xi)$ is a torsion point on $A_\xi  \}$ be infinite, one must in general insist that $dim(X) \geq g$. So, the counterexample does not extend to extensions by $\G_m$ of higher dimensional abelian schemes {\it over curves}.

\section{... in support of Pink's general conjecture.}

 In \cite{Pk}, R. Pink mentions the similarity of Conjecture 6.2 with Y. Andr\'e's result on special points on elliptic pencils \cite{An}, III, p. 9.  Viewing the scheme $B/X$ above   as a ``semi-abelian pencil" over the fixed elliptic curve $E_0$, one could define its special points as the torsion points lying on a fiber which is an isotrivial extension. As mentioned in passing during the proof of the lifting property, the curve $Y = \beta(X)$ even contains infinitely many special points in this sense.
 
 \medskip
 Going further in this direction, we will now construct a mixed Shimura   variety $S(\varphi)$ into which the  image $Y = \beta(X)$ of the Ribet section $\beta$  can be mapped in a natural way. Denoting this map by $i  : Y \rightarrow S(\varphi)$, we have under the hypotheses of \S 1 (or more generally, of Footnote $^{(2)}$ below) : 
 
 \begin{Theorem} The algebraic subvariety $Z = i (Y)$ of the mixed Shimura variety $S(\varphi)$  passes through a Zariski-dense set of special points of $S(\varphi)$, and is indeed a special subvariety of $S(\varphi)$.
 \end{Theorem}
 This, of course, is in full concordance with the prediction of  the general Conjecture 1.3 of \cite{Pk} (more specifically, of the case $d = 0$  of Conjecture 1.1).
 
 \medskip
Here, I will merely give a set-theoretic description of the
construction of $S(\varphi)$.  We fix an integer $g \geq 1$ and  a
totally imaginary quadratic integer $\alpha = - \overline \alpha$, and
denote by $S_0$ (a component of)  the pure Shimura variety
parametrizing abelian varieties $A$ endowed with a  principal
polarization $\psi  : \hat A \rightarrow   A$ (in particular, $\psi =
\hat \psi$), with some level structure, and with an embedding  $j :
\Z[\alpha]  \rightarrow \End(A)$ 
such that $\psi \circ \widehat {j(\alpha)} \circ \psi^{-1} = j(\overline \alpha)$. Let 
 $$\big( {\cal A}, \; \psi, \;  \varphi :=   j(\alpha) \circ \psi \big)$$
  be the corresponding universal  abelian scheme over $S_0$ (in particular, $\varphi = - \hat \varphi : \hat {\cal A} \rightarrow {\cal A}$ is antisymmetric). Then, $S_1 := {\cal A}Ê\times_{S_0} \hat{\cal A}$ is a mixed Shimura variety parametrizing  one-motives of the shape $M : \Z \rightarrow A \times \hat A$, i.e. couples of points $(p, q) \in A \times \hat A$, with $\{A\}Ê\in S_0$, and we can view
  $$S_1(\varphi) = \{ (A, p,q) \in S_1, p = 2\varphi(q) \}$$
as a mixed Shimura subvariety of $S_1$.  Finally, consider the  Poincar\'e bi-extension 
$$ \varpi = (\varpi_1, \varpi_2) : {\cal P}^* \rightarrow  {\cal A}Ê\times_{S_0} \hat{\cal A}.$$
 This is a $\G_m$-torsor  over $S_1 =  {\cal A}Ê\times_{S_0} \hat{\cal A}$, which can again be viewed as a mixed Shimura variety, now parametrizing one-motives $M : \Z \rightarrow B$ of constant and toric ranks equal to 1, where on denoting by $\big(A,p := \varpi_1(M),q := \varpi_2(M)\big)$ the point defined by $\varpi(M)$ in $S_1$,  $B = B_q$ is the extension  of $A$ by $\G_m$ attached to $q$, while the image $b$ of $1 \in \Z$ is a point on $B$ projecting to $p \in A$. Breaking the symmetry between $\cal A$ and $\hat {\cal A}$, we can alternatively consider  $\hat {\cal A} \simeq  \underline{Ext}_{S_0}^1({\cal A}, \G_m)$ as a mixed Shimura variety, and view $\varpi_2 : {\cal P}^* \rightarrow \hat {\cal A}$  as the ``universal" extension ${\cal B}$ of $\cal A$ by $\G_m$, over its parameter space $\hat {\cal A}$. We at last define 
 $$S(\varphi) := \varpi^{-1}(S_1(\varphi))$$
as the mixed Shimura subvariety of ${\cal P}^* \simeq {\cal B}$ whose points parametrize one-motives $M : \Z \rightarrow B$  such that $\varpi_1(M) = 2 \varphi(\varpi_2(M))$.   

\medskip
We now consider an abelian scheme $A/X$ of the type parametrized by $S_0$, over some irreducible algebraic variety $X/\C$. There then exists a unique morphism $i_0 : X \rightarrow S_0$ such that $A/X$ is the pull-back of ${\cal A}$ under $i_0$. We fix a section $q : X \rightarrow \hat A$, corresponding to a semi-abelian scheme $\pi : B \rightarrow A$ over $X$, and perform Ribet's construction\footnote {~ When $g = 1$ as in the first Section, $S_0$ is reduced to a CM point $\{E_0\}$. Taking into account Remark 1.(ii), we are now proving Theorem 2 for any $ g \geq 1$. However, in the case $g > 1$, we must add the following   hypotheses to its statement : the base $X/\C$ is an irreducible variety of dimension $\geq g$, and

(a)  the image $i_0(X)$ meets  the set of CM  points of $S_0$ Zariski-densely.  

(b)   $q(X)$ meets the set of torsion points of the various CM fibers of $\hat A$ Zariski-densely;

\noindent
For the sake of simplicity, we will  {\it assume}, in what follows, {\it that $i_0(X)$ is a Shimura subvariety of $S_0$} (more or less $\Leftrightarrow$ (a) under Andr\'e-Oort), and {\it that   $q : X \rightarrow \hat A$ dominates the ``generic case" $id : \hat A \rightarrow \hat A$} ($\Rightarrow$ (b)). Furthermore, the latter hypothesis implies the useful (but not necessary) property that :

(c)  $q$ factors through no    proper closed subgroup scheme of $ \hat A/X$.}, yielding  sections $p = 2 \varphi \circ q : X \rightarrow A,  \; \beta : X \rightarrow B$, with $p = \pi \circ \beta$.   By the universal property of $S(\varphi)$, there exists a unique morphism  
$$i : X \rightarrow S(\varphi)$$
 above $i_0$ such that the smooth  $X$-one-motive $M : \Z \rightarrow B$ defined by $\beta$ is the pull-back under $i$ of the universal one-motive ${\cal M} : \Z \rightarrow {\cal B}$. We again denote by $i : B \rightarrow S(\varphi) \subset {\cal B}$ the extension of $i$ to $B/X$. The image $Z := i(Y) \subset S(\varphi)$ of $Y = \beta(X) \subset B$ is the algebraic subvariety of  $S(\varphi)$ to be studied for Theorem 2.
 
 \medskip
 We first check that $Z$ contains a Zariski-dense set of special points of $S(\varphi)$. By definition, these points represent complex one-motives such that the underlying abelian variety is CM, and whose  weight filtrations  are totally split up to isogeny.  By the first hypothesis made in Footnote $^{(2)}$, the projection $i_0(X)$ of $Z$ to $S_0$ passes through a Zariski-dense set of CM points. In the fiber of any such  point, the projection of $Z$ to $S_1(\varphi)$ passes through a Zariski-dense set of torsion points $(p = 2 \varphi(q), q)$, because of  the second hypothesis made in this foot-note. The lifting property established in \S 1 (more accurately, the sharper version mentioned in passing) now shows that $Z$ does meet the set of special points of $S(\varphi)$ Zariski-densely.
 
 \medskip
 It remains to show that $Z$ is a special subvariety of $S(\varphi)$. We will check this Hodge-theoretically. Denote by $G$, resp. $P$,  the generic Mumford-Tate group of the Shimura  variety $i_0(X)$, resp. of the mixed Shimura subvariety of $S(\varphi)$ lying above $i_0(X)$.  Then, $P$ is the semi-direct product of  its unipotent radical $W_{-1}(P)$ by $G$, and $W_{-1}(P)$ is the semi-direct product of  $W_{-2}(P) \subset \G_a$ by the  vectorial group $Gr_{-1}(P) \subset \G_a^{2g}$  : notice that $Gr_{-1}(P)$ is the unipotent radical of the generic Mumford-Tate group  of $S_1(\varphi)$, and the inclusion  $Gr_{-1}(P) \subset \G_a^{2g}$ follows from the linear dependence relation $p = 2 \varphi(q)$.  Similarly, let  $P_z \subset P$ denote the  Mumford-Tate group of  a sufficiently general point $z$ in $Z$; then, $P_z$  is the semi-direct product  of  its unipotent radical $W_{-1}(P_z)$ by $G$, and $W_{-1}(P_z)$ is the semi-direct product of $W_{-2}(P_z) \subset  W_{-2}(P)$ by the vectorial group $Gr_{-1}(P_z) \subset Gr_{-1}(P)$.
 
 \smallskip
 Now,  by the second hypothesis of the footnote, $\Z. q(X)$ is Zariski-dense in $\hat A$, and Proposition 1 of \cite{Ay} shows that $Gr_{-1}(P_z) = Gr_{-1}(P) = \G_a^{2g}$.  On the other hand, Theorem 1 of \cite{Bd} shows that $W_{-2}(P) = \G_a$, while $W_{-2}(P_z) = \{0\}$, and more precisely, that {\it for any point $s \in S(\varphi)(\C)$ whose Mumford-Tate group $P_s$ satisfies  $Gr_{-1}(P_s) = \G_a^{2g}$, we have}
 $$W_{-2}(P_s) = \{0\} \Leftrightarrow \exists \, \tilde s \in [s], \, \tilde s \in Z,$$ 
where $[s]$ denotes the Hecke orbit of $s$. In other words, up to isogenies, the  points  of $Z$ are characterized by the existence of an exceptional Hodge tensor in their Betti realization, which does not exist at the generic point of $S(\varphi)$. So, $Z$ is indeed a   special subvariety  of $S(\varphi)$.

 \bigskip
 \noindent
 {\bf Remark 2} : (i) This section shows that the special subvarieties of the mixed Shimura variety $\cal B$ 
do not necessarily  correspond to families of semi-abelian subvarieties, so that Theorem 5.7 and 6.3 of \cite{Pk} must be modified. Similarly, the mixed Shimura varieties of Hodge type $\{W_{-2}(P) = 0\}$ are not necessarily Kuga fiber varieties (compare \cite{Mi}, Example 1.10, and \cite{Bd}, Remark (v)).

\smallskip
\qquad (ii) The strange divisibility properties of Ribet points alluded to in the introduction of \S 1 are precisely reflected by the  $\ell$-adic analogue of the vanishing of $W_{-2}(P_z), z \in S(\varphi)(\overline \Q)$; cf. \cite{Bd}, Theorem 1.(ii). Actually, K. Ribet gives in \cite{Ri}  an explicit  description of  the  exceptional  Galois invariant tensor occuring in the $\ell$-adic cohomology of these one-motives, which applies to all their realizations.  
 
 \smallskip
\qquad  iii)  We refer to  \cite{Bn} for a complete description of the unipotent radical $W_{-1}(P)$ of the Mumford-Tate group of one-motives of higher toric or constant ranks.  When both ranks  are equal to 1, the fact that  $W_{-2}(P)$ vanishes in the case  of antisymmetrically self-dual one-motives is an exercise in goup theory, cf. \cite{Bp},  Lemme 6.

\section{Appendix, by Bas Edixhoven}

We go back to the setting and the notations $E_0, X, q, \varphi =
-\hat \varphi, \pi : B \rightarrow E$ of \S\,1, but to make later
comparisons easier, we henceforth consider the extension
$$B = B_{2q}$$ 
of $E$ by $\G_m$ given by $2q\in \hat{E}(X)$, and
denote by $\beta_R$ in $B_{2q}(X)$ the Ribet section which the
construction of \S\, 1 attaches to the data $\{2q, \varphi \}$; in
particular, the projection $p := \pi \circ \beta_R = \varphi(2q) -
\hat \varphi(2q) = 2 \varphi(2q))$ of $\beta_R$ on $E(X)$ is now twice
the section considered in $\S 1$.

\medskip
Interpreting $B_{2q}/X$ as the generalized jacobian $\Pic^0_{C/X}$ of
a singular curve $C/X$ with normalization $\hat{E} = \Pic^0_{E/X} =
\hat E_0 \times X/X$, we will here construct\footnote{~I thank Lenny
  Taelman for the suggestion to pay more attention to the symmetry
  between the points $q$ and $-q$ to be identified and the divisor
  $\beta_a$ that gives the section~$\beta_J$.}  a concrete section
$\beta_J$ (with ``J'' for ``Jacobian'') of $B_{2q}/X$, which enjoys
all the properties\footnote{~ In particular, $\beta_J$ does not factor
  through any proper closed subgroup scheme of $B$.  For a
  discussion relating $\beta_J$ and $\beta_R$, see Remark 3.ii below.}
of the Ribet section $\beta_R$, and thereby provides a new proof of
Theorem~1, but for which the ``lifting property" takes the following
sharper form \;:

\begin{Theorem} 
For any $\xi \in X$, the image of $\beta_J(\xi)$ under $\pi : B_\xi
\rightarrow E_\xi \simeq E_0$ satisfies $\pi(\beta_J(\xi)) =
\pi(\beta_R(\xi)) :=  \,p(\xi) $. And if $p(\xi)$ is a torsion
point of $ E_0$ of order $n$, with $n$ prime to $2\, \deg(\varphi)
\deg(\varphi+\psi)\deg(\varphi-\psi)$, then $\beta_J(\xi)$ is a
torsion point of $B_\xi$, of order dividing~$n^2$.
\end{Theorem}
 
\medskip
Here, $\psi = \hat \psi: \hat E \rightarrow E$ denotes the standard
principal polarization: its inverse sends a point $P$ to the class of
the divisor $(P)-(0)$. Without loss of generality, we assume that we
are in the generic case $X = \hat E_0, q = \mathrm{id}$, with $K :=
\C(X) = \C(\hat E_0)$, but we keep to the notation $X$ to indicate
that some points of $X$ may have to be removed in the constructions
which follow.  We denote by $Q = \psi \circ q : X \rightarrow E$ the
section of $E/X$ such that $q \in \hat E(X)=\Pic^0(E/X)/\Pic(X)$ is
represented by the divisor $(Q) -(0)$ on~$E$. By biduality, the
section $Q$ in $E(X)=\Pic^0(\hat E/X)/\Pic(X)$ is represented by the
divisor $(q)-(0)$ on~$\hat E$.

\medskip
Let $C/X$ be the singular curve over $X$ obtained by identifying the
disjoint sections $q$ and $-q$ of~$\hat E$ (we remove $\hat{E}_0[2]$
from~$X$). As a set, it is the quotient of $\hat E$ by the equivalence
relation generated by $q(\xi)\sim -q(\xi)$ with $\xi$ ranging
over~$X$. The topology on $C$ is the finest one for which the quotient
map $\quot\colon\hat E\to C$ is continuous: a subset $U$ of $C$ is
open if and only if $\quot^{-1}U$ is open in~$\hat E$. The regular
functions on an open set $U$ of $C$ are the regular functions $f$ on
$\quot^{-1}U$ such that $f(q(\xi))=f(-q(\xi))$ whenever
$\quot(q(\xi))$ is in~$U$. It is proved in Thm.~5.4 of \cite{Ferrand}
that this topological space with sheaf of $\C$-valued functions is
indeed an algebraic variety over~$\C$.  In categorical terms,
$\quot\colon \hat E\to C$ is the equalizer of the pair of morphisms
$(q,-q)$ from $X$ to~$\hat E$.

\smallskip
The curve $C\to X$ is a family of singular curves, each with an
ordinary double point; it is semi-stable of genus two
(see~\cite[9.2/6]{BLR}). Its normalization is $\quot\colon \hat E \to
C$. Its generalized jacobian $B:=\Pic^0_{C/X}$ is described in
\cite{BLR}, 8.1/4, 8.2/7, 9.2/1, 9.3/1. As $C\to X$ has a section (for
example $\bar{0}:=\quot\circ 0$ and $\bar{q}:=\quot\circ q$), we have,
for every $T\to X$, that $B(T)=\Pic^0(T\times_XC/T)/\Pic(T)$, where
$\Pic^0(T\times_XC/T)$ is the group of isomorphism classes of line
bundles on $T\times_XC$ that have degree zero on the fibres of
$T\times_XC\to T$. The group $\Pic(T)$ is contained as direct summand
in $\Pic^0(T\times_XC/T)$ via pullback by the projection
$T\times_XC\to X$ and a chosen section. In particular, a divisor $D$
on $C$ that is finite over~$X$, disjoint from $\bar{q}(X)$ and of
degree zero after restriction to the fibers of $C\to X$ gives the
invertible $\calO_C$-module $\calO_C(D)$ that has degree zero on the
fibers and therefore gives an element denoted $[D]$ in $B(X)$.

\smallskip
For $\xi$ in $X$, the fiber $B_\xi$ is, as abelian group, the group
$\Pic^0(C_\xi)$. In terms of divisors this is the quotient of the
group $\Div^0(C_\xi)$ of degree zero divisors with support outside
$\{\bar{q}(\xi)\}$ by the subgroup of principal divisors $\divisor(f)$ for
nonzero rational functions $f$ in $\C(C_\xi)^\times$ that are regular
at~$\bar{q}(\xi)$. As $C_\xi-\{\bar{q}(\xi)\}$ is the same as $\hat
E_0-\{q(\xi),-q(\xi)\}$, $\Div^0(C_\xi)$ is the group of degree zero divisors
on $\hat E_0$ with support outside $\{q(\xi),-q(\xi)\}$. An element $f$ of
$\C(C_\xi)^\times$ that is regular at $\bar{q}(\xi)$ is an element of
$\C(\hat E_0)^\times$ that is regular at $q(\xi)$ and $-q(\xi)$ and satisfies
$f(q(\xi))=f(-q(\xi))$. This gives us a useful description of~$B_\xi$.

\medskip
The normalization map $\quot\colon \hat E \rightarrow C$
induces a morphism of group schemes over $X$
\[
\pi\colon  B = \Pic^0_{C/X} \rightarrow \Pic^0_{\hat E/X} =  E/X,
\]
and identifies $B$ with the extension of $E$ by $\G_m$ given by the
section 
$2q \in \hat{E}(X)$. For $\xi$ in $X$, the class $[\delta]$ in
$B_\xi$ of a divisor $\delta \in \Div^0(C_\xi)$ lies in the kernel
$\C^\times$ of $\pi_\xi$ if and only if there exists $f \in \C(\hat
E_0)^\times$ such that $\delta = \divisor(f)$ on $\hat E_0$, and it is
then a torsion point in $\C^\times$ if and only if the quotient
$f(q(\xi))/f(-q(\xi)) \in \C^\times$, which does not depend on the
choice of~$f$, is a root of unity.

\smallskip
We recall that for an elliptic curve $\cal E$ over an algebraically closed
field $k$, $\psi\colon \hat{\cal E}\to \cal E$ the standard polarization and $u$ in
$\End(\cal E)$, the pullback map $u^*$ on $\Div(\cal E)$ induces $\hat{u}$ in
$\End(\hat{\cal E})$, the dual of~$u$, and then
$\bar{u}:=\psi\hat{u}\psi^{-1}$ in $\End(\cal E)$ is called the Rosati-dual
of~$u$; it is characterized by the property that in $\End(\cal E)$ we have
$\bar{u}u=\deg(u)$ and $u+\bar{u}\in\Z$. Also, the pushforward map
$u_*$ on $\Div(\cal E)$ induces an element still denoted $u_*$ in
$\End(\hat{\cal E})$ such that $\psi^{-1}u=u_*\psi^{-1}$ and $u_*u^*$ is
multiplication by $\deg(u)$ in $\End(\hat{\cal E})$. Hence $u_*$ is the
Rosati dual of~$u^*$. We have $\hat{\hat{u}}=u$, $\bar{\bar{u}}=u$ and
$\bar{u}_*=u^*$ in $\End(\hat{\cal E})$. For $f$ a nonzero rational
function on $\cal E$ and $u\neq 0$ we have $u^*\divisor(f)=\divisor(f\circ
u)$, and $u_*\divisor(f)=\divisor(\Norm_u(f))$, where $\Norm_u\colon
k({\cal E})^\times\to k({\cal E})^\times$ is the norm map along~$u$. Of course, all
this applies to $E_0$ and $\hat{E}_0$ over $\C$, and to $E$ and $\hat{E}$ over
the algebraic closure of the function field $K$ of~$X$.

\smallskip
We will use Weil reciprocity: for $f$ and $g$ nonzero rational
functions on a nonsingular irreducible projective curve $\cal E$ over an
algebraically closed field $k$ such that $\divisor(f)$ and $\divisor(g)$
have disjoint supports, one has $f(\divisor(g))=g(\divisor(f))$, 
where for $D=\sum_PD(P){\cdot}P$ a divisor on $\cal E$ one defines
$f(D)=\prod_Pf(P)^{D(P)}$, cf. \cite{Se}, III, Prop. 7. 
In Remark 3.(i) after the proof, we will also use the Weil pairing
on~$\cal E$. For $n$ a positive integer and $P$ and $Q$ in $\Pic^0({\cal E})[n]$
the element $e_n(P,Q)$ in $\mu_n(k)$ is defined as follows. Let 
$D_P$ and $D_Q$  in $\Div^0(\cal E)$ be disjoint divisors representing $P$
and~$Q$. Let $f$ and $g$ be in $k({\cal E})^\times$ such that
$nD_P=\divisor(f)$ and $nD_Q=\divisor(g)$. Then $e_n(P,Q)=f(D_Q)/g(D_P)$. For
$n$ invertible in $k$ this pairing $e_n$ is a perfect alternating
pairing.

\medskip
To define $\beta_J = \beta_J(a)$, we let $a\in\End(\hat{E}_0)$ be
\emph{any} endomorphism such that $a^3-a\neq 0$. We set
$\alpha:=\hat{a}$ in $\End(E_0)$ and
$\varphi:=\alpha\circ\psi\colon\hat{E_0}\to E_0$.  Just as in
\cite{JR}, we will not need that $\hat{\varphi}=-\varphi$.  However,
if $\hat{\varphi}=\varphi$ then $\pi\beta_J$ will be zero in
$\hat{E}(X)$, so we will insist that $\hat{\varphi}-\varphi\neq 0$.
For $s$ in $\hat{E}(X)$ we let $(s)$ denote the relative divisor that
it gives on $\hat{E}$, and, if $(s)$ is disjoint from the singular
locus $\bar{q}(X)$ of $C$, also the relative divisor on $C$ that it
gives.

\smallskip
We set:
\[
\beta_a := a^*\big((q)-(-q)\big) - a_*\big((q)-(-q)\big) 
\quad\text{in $\Div^0(\hat{E}/X)$}. \qquad \qquad (1)
\]
The reader should note the antisymmetric use of both $a^*$ and $a_*$ in
the definition of~$\beta_a$. To get its support disjoint from
$\bar{q}(X)$ we remove from $X$ the finite set $\ker(a^2-1)$: note
that, for $\xi\in X$,  $\beta_a(\xi)$ and $(q(\xi))-(-q(\xi))$ are not
disjoint if and only if $aq(\xi)=q(\xi)$ or $aq(\xi)=-q(\xi)$. We
can now also view $\beta_a$ as element of $\Div^0(C/X)$, and we set:
\[
\beta_J := [\beta_a] \quad\text{in $B(X)$.}
\]
The image of $\pi\beta_J$ of $\beta_J$ in $E(X)$ is the class of
the divisor $\beta_a$ on $\hat{E}$, hence we have, 
on denoting by $\simeq$ linear equivalence on $\Div^0(\hat E/X)$ :
\begin{align*}
\pi\beta_J & = \bar{a}_*\big((q)-(-q)\big)-a_*\big((q)-(-q)\big)  \\
& \simeq  (\bar{a}q)-(\bar{a}(-q)) - ((aq)-(a(-q)) \\
& \simeq \left(2(\bar{a}-a)q\right)-(0) 
  = 2(\alpha-\bar{\alpha})Q \quad\text{in $E(X)$.}
\end{align*}
In particular, since $\bar{a}-a$ is nonzero, then $\pi\beta_J$ is not
torsion in~$E(X)$, and in fact
$$ 
\pi \beta_J = \big((2\bar a q) - (0)\big) - \big((2aq) -(0)\big) =
\varphi(2q) - \hat \varphi(2q) = p = \pi \beta_R,
$$
as was to be checked for the first part of Theorem 3.

\medskip
Now we start the proof of the second part of Theorem~3.  Let $n$ be a
positive integer prime to $2 \deg(\varphi) \deg(\varphi + \psi)
\deg(\varphi - \psi) = 2 \deg(a(a^2 -1))$, and let $\xi\in X$ be a
point such that $p(\xi)\in E_\xi= E_0$ is a torsion point of order
$n$. Actually, we will {\it assume that $q(\xi)$ is a point of order
  $n$ in $\hat{E}_\xi=\hat{E}_0$} : in the antisymmetric case
considered in Theorem 3, we have $p(\xi) = 4 \varphi(q(\xi))$, and the
two conditions are equivalent by the primality assumption. To ease
notations, we now drop the mention of $\xi$, writing $E, B, q,\ldots$
instead of $E_\xi, B_\xi, q(\xi),\ldots$.

\medskip
As $nq=0$ in $\hat{E}_0$, we have $n\pi\beta_J=0$ in~$E_0$. This means
that $n\beta_a$ is a principal divisor on~$\hat{E_0}$. Let
$f\in\C(\hat{E}_0)^\times$ be such that $\divisor(f)=n(q)-n(-q)$ in
$\Div(\hat{E}_0)$. Then we have, on $\hat{E}_0$:
\begin{align*}
\divisor(f\circ a) & = a^*\divisor(f) = a^*\left(n(q)-n(-q)\right),\\
\divisor(\Norm_a(f)) & = a_*\divisor(f) = a_*\left(n(q)-n(-q)\right).
\end{align*}
We define:
\[
g_a:=(f\circ a)/\Norm_a(f)\quad\text{in $\C(\hat{E}_0)^\times$.}
\]
Then we have:
\[
n\beta_a = \divisor(f\circ a)-\divisor(\Norm_a(f)) = \divisor(g_a) 
\quad\text{on $\hat{E}_0$.} 
\]
This means that $n[\beta_J]$ in $B$ is the element $g_a(q)/g_a(-q)$
of~$\C^\times$. As the divisor of $f$ has support disjoint from that
of $g_a$ and of $a^*\divisor(f)$ and $a_*\divisor(f)$, Weil
reciprocity gives us:
\begin{align*}
\left(\frac{g_a(q)}{g_a(-q)}\right)^n & = g_a(\divisor(f)) =
f(\divisor(g_a)) = f(\divisor(f\circ a)-\divisor(\Norm_a(f)))\\
& = \frac{f(\divisor(f\circ a))}{f(\divisor(\Norm_a(f)))} = 
\frac{(f\circ a)(\divisor(f))}{f(a_*\divisor(f))} = 
\frac{f(a_*\divisor(f))}{f(a_*\divisor(f))} = 1.
\end{align*}
This finishes the proof of Theorem~3.

\bigskip
\noindent
{\bf Remark 3} : i) We have shown that for $\xi$ torsion in $X \subset
\hat{E}_0$ and $n$ its order, we have $n^2\beta_J(\xi)=0$
in~$B_\xi$. We will show that for $n>1$ odd and prime to
$\deg((a^2-1)(\bar{a}-a))$ there exist $\xi$ in $X$ with $p(\xi)$ of
order $n$ in $E_0$ such that the order of $\beta_J(\xi)$ equals~$n^2$.

Let $n$ be such an integer. Recall the notation $Q = \psi(q)$, and let
$\xi$ in $\hat{E}_0$ be of order $n$ such that
$e_n(2(\alpha-\bar{\alpha})Q(\xi),2Q(\xi))$ is of order $n$
in~$\C^\times$. Such $\xi$ exist because $\alpha-\bar{\alpha}$ is an
automorphism of $E_0[n]$ that is not scalar multiplication by an
element of $\Z/n\Z$. And such a $\xi$ is in $X$ because $n$ is prime
to $2\deg(a^2-1)$. By construction, $p(\xi) =
2(\alpha-\bar{\alpha})Q(\xi)$ is represented by the divisor
$\beta_a(\xi)$, and $2Q(\xi)$ is represented by
$(q(\xi))-(-q(\xi))$. Again, let us drop the $\xi$'s from our
notation. Then $n\beta_a=\divisor(g_a)$, and
$n((q)-(-q))=\divisor(f)$. We first compute:
\begin{align*}
f(\beta_a) & = f(a^*((q)-(-q))-a_*((q)-(-q))) = 
\frac{(\Norm_af)((q)-(-q))}{(f\circ a)((q)-(-q))} \\
& = g_a^{-1}((q)-(-q)) = \frac{g_a(-q)}{g_a(q)}.
\end{align*}
Hence:
\[
e_n(2(\alpha-\bar{\alpha})Q,2Q) = \frac{g_a((q)-(-q))}{f(\beta_a)} = 
\frac{g_a(q)}{g_a(-q)f(\beta_a)} =  \left(\frac{g_a(q)}{g_a(-q)}\right)^2.
\]

\smallskip
ii) The Ribet section $\beta_R$ and the present section $\beta_J$ of
$B/X$ differ by an element of $\G_m(X) $, and it is natural to ask
whether they are actually equal.  In this respect, notice that by
Theorem 3, $\beta_J(X)$ passes through a Zariski-dense set of special
points of the mixed Shimura variety $S(\varphi)$ of \S 2. If we assume
Pink's general Conjecture 1.3 of \cite{Pk}, $\beta_J(X)$ must then be
a component of the Hecke orbit of the special subvariety $Z$ defined
by $\beta_R(X)$. Hence, at least conjecturally, the section $ \beta_R
- \beta_J $ of $B/X$ is a root of unity.

 \bigskip
 \noindent
{\it  Authors' addresses} : 
 
 \smallskip
 D.B. : Institut de Math\'ematiques de Jussieu, bertrand@math.jussieu.fr
 
 B.E. : Mathematisch Instituut, Univ. Leiden, edix@math.leidenuniv.nl
 
 \bigskip
 \noindent 
{\it  AMS Classification} : 11G18, 14H40, 14K12 .

\noindent
 {\it Key words}  : semi-abelian varieties, Manin-Mumford conjecture, one-motives, generalized jacobians, special points, mixed Shimura varieties.

\end{document}